\documentclass{amsart}
\usepackage{amssymb,latexsym}

\newtheorem*{named}{Main Theorem}

\newtheorem{theorem}{Theorem}
\newcommand{\bt}{\begin{theorem}}
\newcommand{\et}{\end{theorem}}

\newtheorem{lemma}{Lemma}
\newcommand{\bl}{\begin{lemma}}
\newcommand{\el}{\end{lemma}}

\newtheorem{corollary}{Corollary}
\newcommand{\bc}{\begin{corollary}}
\newcommand{\ec}{\end{corollary}}
\newcommand{\Z}{\ensuremath{\mathbf Z}}
\newcommand{\beq}{\begin{equation}}
\newcommand{\eeq}{\end{equation}}
\newcommand{\benum}{\begin{enumerate}}
\newcommand{\eenum}{\end{enumerate}}
\DeclareMathOperator{\card}{\text{card}}
\DeclareMathOperator{\ord}{\text{ord}}
\newtheorem{problem}{Problem} 
\newcommand{\bprob}{\begin{problem}}
\newcommand{\eprob}{\end{problem}}

\begin{document}

\title[Bi-Lipschitz metrics and additive number theory]{Bi-Lipschitz equivalent metrics on groups, and a problem in additive number theory}
\author{Melvyn B. Nathanson}
\address{CUNY (Lehman College and the Graduate Center)}
\email{melvyn.nathanson@lehman.cuny.edu}

\date{\today}



\thanks{This paper was supported in part by a PSC-CUNY Research Award, and was written while M.B.N. was a visiting fellow at  Princeton University.  He thanks the Princeton math department for its hospitality.}

\begin{abstract}
There is a standard ``word length'' metric canonically associated to any set of generators for a group.  In particular, for any integers $a$ and $b$ greater than 1,  the additive group \Z\ has generating sets 
$\{ a^i \}_{i=0}^{\infty}$ and $\{b^j\}_{j=0}^{\infty}$ with associated metrics $d_A$ and $d_B$, respectively.  It is proved that these metrics are bi-Lipschitz equivalent if and only if there exist positive integers $m$ and $n$ such that $a^m = b^n$.
\end{abstract}

\maketitle

\section{Groups as metric spaces}
Two metrics $d_A$ and $d_B$ on the same set $X$ are called \emph{bi-Lipschitz equivalent metrics} if there exists a number $K \geq 1$ such that 
\beq  \label{Bi-Lipschitz:metric}
\frac{1}{K}d_A(x_1,x_2) \leq d_B(x_1,x_2) \leq K d_A(x_1,x_2)
\eeq
for all $x_1,x_2 \in X$.  In this paper we study metrics on the additive group \Z\ of integers that are defined in a natural way by geometric progressions, and give a necessary and sufficient condition for two such metrics to be bi-Lipschitz equivalent.

More generally, two metric spaces $(X, d_X)$ and $(Y,d_Y)$ are called \emph{bi-Lipschitz equivalent} if there exists a function $f$ from $X$ onto $Y$ such that
\[
\frac{1}{K}d_X(x_1,x_2) \leq d_Y( f(x_1),f(x_2) ) \leq K d_X(x_1,x_2))
\]
for some number $K \geq 1$ and all $x_1,x_2 \in X$.  The function $f$ is called a \emph{bi-Lipschitz equivalence}; it is necessarily a homeomorphism.   In particular, the metrics $d_A$ and $d_B$ on a set $X$ are bi-Lipschitz equivalent metrics if and only if the identity map from $(X,d_A)$ to $(X, d_B)$ is a bi-Lipschitz equivalence.

Let $G$ be a group with identity $e$, and let $A$ be a set of generators for $G$.  
The set $A$ may be finite or infinite.  Let $A^{-1} = \{ a^{-1}: a \in A \}$.   A  \emph{word with respect to $A$} is a finite product of the form $a_1a_2 \cdots a_n$, where $a_i \in A \cup A^{-1}$ for $i=1,\ldots, n$.  We call $n$ the \emph{length} of this word.  An element $x \in G\setminus \{ e \}$ has \emph{length} $\ell_A(x) = n$ with respect to the generating set $A$ if $n$ is the least positive integer such that $x$ can be written as a word with respect to $A$ of length $n$.  We define $\ell_A(e) = 0$.  Every element of $G$ has finite length because $A$ generates $G$.  

The length function has the following properties.  
First, $\ell_A(x) = 0$ if and only if $x=e$.  Second, $\ell_A(x) = \ell_A(x^{-1})$ for all $x \in G$.  Third, there is the subadditivity condition
\beq  \label{quasi:LengthIneq}
\ell_A(xy) \leq \ell_A(x) + \ell_A(y)
\eeq
for all $x,y \in G$.  

We use the length function associated with the generating set $A$ to construct a distance function $d_A$ on the group $G$.   For all $x, y \in G$, let $d_A(x,y)$ be the length of the group element $xy^{-1}$, that is,
\[
d_A(x,y) = \ell_A(xy^{-1})
\]  
for all $x,y\in X$.  
In particular, $d_A(x,e) = \ell_A(x)$ for all $x \in G$.   It follows that  $d_A(x,y) = 0$ if and only if $\ell_A(xy^{-1}) = 0$, that is, if and only if $xy^{-1} = e$ or $x=y$.  Similarly, 
\[
d_A(x,y) = \ell_A(xy^{-1}) = \ell\left( (xy^{-1})^{-1} \right) 
= \ell\left( yx^{-1}\right) =  d_A(y,x).
\]
By inequality~\eqref{quasi:LengthIneq}, we have
\begin{align*}
d_A(x,z) & = \ell\left( xz^{-1} \right) = \ell\left( xy^{-1}yz^{-1} \right) \\
& \leq   \ell\left( xy^{-1} \right)+  \ell\left( yz^{-1} \right) \\
& =  d_A(x,y) + d_A(y,z)
\end{align*}
for all $x,y,z \in G$, and so $d$ satisfies the triangle inequality.  
Thus, $(G,d_A)$  is a metric space.  We call  $d_A$ the \emph{metric associated to the generating set $A$}.

\bl    \label{quasi:lemma:MetricEquiv}
Let $A = \{a_i\}_{i\in I}$ and $B = \{ b_j\}_{j\in J}$ be generating sets for the group $G$.  Let $\ell_A$ and $\ell_B$ be the corresponding length functions, and $d_A$ and $d_B$ the associated  metrics on $G$.  The following conditions are equivalent:
\benum
\item[(i)]
$\sup\{ \ell_A(b_j) : j \in J \} < \infty$
and
$\sup\{ \ell_B(a_i) : i \in I \} < \infty$
\item[(ii)]
the metrics $d_A$ and $d_B$ are bi-Lipschitz equivalent
\eenum
\el

\begin{proof}
If $K_B = \sup\{ \ell_A(b_j) : j\in J \} < \infty$, then every generator $b_j \in B$ can be represented as product of at most $K_B$ elements of $A \cup A^{-1}$.  Taking the inverse of this representation, we see that the inverse generator $b_j^{-1}$ can also be represented as a product of at most $K_B$ elements of $A \cup A^{-1}$.  Let $x,y \in G$.  If $d_B(x,y) = n$, then $xy^{-1}$ is a product of $n$ elements of $B \cup B^{-1}$.  Writing each of these as a product of at most $K_B$ elements of $A \cup A^{-1}$, we obtain a representation of $xy^{-1}$ as the product of at most $K_B n$  elements of $A \cup A^{-1}$, and so $d_A(x,y) \leq K_B n = K_B d_B(x,y)$.  Similarly, if $K_A = \sup\{ \ell_B(a_i) : i \in I \} < \infty$, then every element in $A \cup A^{-1}$ can be represented as a product of at most $K_A$ elements of $B \cup B^{-1}$, and $d_B(x,y) \leq K_A d_A(x,y)$.  Thus, if $K_A<\infty$ and $K_B < \infty$, then inequality~\eqref{Bi-Lipschitz:metric} holds with $K = \max(K_A,K_B)$, and the metrics $d_A$ and $d_B$ are equivalent.
This proves that (i) implies (ii).

Conversely, if the metrics are $d_A$ and $d_B$ are bi-Lipschitz equivalent, then there exists a number $K \geq 1$ such that inequality~\eqref{Bi-Lipschitz:metric} holds.  For every generator $a_i \in A$ we have
\[
\ell_B(a_i) = d_B(a_i,e) \leq K d_A(a_i,e) = K\ell_A(a_i) = K  
\]
and so $\sup\{ \ell_B(a_i) : i \in I \} \leq K < \infty$.  Similarly, $\sup\{ \ell_A(b_j) : j \in J \} \leq K < \infty$.   Therefore, (ii) implies (i).  This completes the proof.
\end{proof}

\bc
Let $A = \{a_i\}_{i\in I}$ and $B = \{ b_j\}_{j\in J}$ be generating sets for the group $G$.  If $A$ and $B$ are finite sets, then the associated  metrics  $d_A$ and $d_B$ are bi-Lipschitz equivalent.
\ec

\begin{proof}
It suffices to observe that the finite sets of numbers  
$\{ \ell_B(a_i) : i \in I \}$ and $\{ \ell_A(b_j) : j \in J \}$  have finite upper bounds.
\end{proof}

\bc
Let $A = \{a_i\}_{i\in I}$ and $B = \{ b_j\}_{j\in J}$ be generating sets for the group $G$.  If $A$ is a finite set and $B$ is an infinite set, then the associated  metrics  $d_A$ and $d_B$ are not bi-Lipschitz equivalent.
\ec

\begin{proof}
If $|A| = r < \infty$, then $|A \cup A^{-1}| \leq 2r$, and for every positive integer $s$ there are less than $(2r)^{s+1}$  words with respect to $A \cup A^{-1}$ of length at most $s$.  Since $B$ is infinite, it follows that there are infinitely many generators $b_j \in B$ with $\ell_A(b_j) > s$, and so $\sup\{ \ell_A(b_j) :  j \in J \} = \infty$.  
Therefore, the metrics $d_A$ and $d_B$ are not bi-Lipschitz equivalent.  
\end{proof}

It remains to determine when the metrics associated with different infinite generating sets for a group are bi-Lipschitz equivalent.   This is an open problem even for \Z, the additive group of integers, for which the generating sets are the sets of relatively prime integers.  We shall prove the following theorem, which determines when the metrics associated with infinite geometric sequences of integers are bi-Lipschitz equivalent.

\begin{named}
Let $a$ and $b$ be integers greater than 1, and consider the additive  group \Z\ with generating sets  $A = \{a^i\}_{i=0}^{\infty}$ and $B = \{b^j\}_{j=0}^{\infty}$.  Let $d_A$ and $d_B$ be the metrics on \Z\ associated with the generating sets $A$ and $B$, respectively.  These metrics are bi-Lipschitz equivalent if and only if there exist positive integers $m$ and $n$ such that $a^m = b^n$.
\end{named}

\section{Representations of integer powers to various integer bases}

Let $a$ be an integer greater than 1.    
Every nonnegative integer $n$ has a unique $a$-adic representation  
\beq    \label{quasi:a-adic}
n = \sum_{i=0}^{\infty} \delta_i a^i 
\eeq
where $\delta_i \in \{ 0,1,2,\ldots, a-1\}$ for all nonnegative integers $i$, and $\delta_i = 0$  for all sufficiently large $i$.  For integers $u < v$ we denote by  $[u,v)$ the interval of integers $\{u,u+1,\ldots, v-1\}$.  An interval $[u,v)$ is called an \emph{$a$-adic block} for the positive integer $n$ if, in the $a$-adic expansion~\eqref{quasi:a-adic}, we have $\delta_i \neq 0$ for all $i \in [u,v)$.    The $a$-adic block $[u,v)$  for the positive integer $n$ is called a \emph{maximal $a$-adic block} if either $u=0$ and $\delta_{v} = 0$, or $u\geq 1$ and $\delta_{u-1} = \delta_v = 0$.  We define the \emph{maximal $a$-adic block function} $M_A(n)$ as the number of maximal $a$-adic blocks in the $a$-adic expansion of $n$.  For example, if $I$ is a set of $k$ nonnegative integers, if $\delta_i \in \{1,2,\ldots, a-1\}$ for $i\in I$, and if $n = \sum_{i\in I} \delta_i a^i$, then 
$
M_A\left(n \right)\leq k.
$
Moreover, $M_A\left(n \right) = k$ if and only if no two elements of $I$ are consecutive.

\bl     \label{quasi:lemma:binary}
Let $a \geq 2$ and $r \geq 1$.  Let $J = \{ j_i\}_{i=0}^r$ be a strictly decreasing sequence  of nonnegative integers and let $\delta_{j_i} \in \{1,2,\ldots, a-1\}$ for $i = 0,1,\ldots, r$.  If 
\[
n_J = \delta_{j_0} a^{j_0} - \sum_{i=1}^r \delta_{j_i} a^{j_i}
\]
then 
\[
M_A(n_J) \leq  r
\]
and
\[
n_J =  \left( \delta_{j_0}  - 1 \right) a^{j_0}
+ \sum_{i=j_r +1}^{j_0 -1} \delta'_{i} a^{i} +   \delta'_{j_r} a^{j_r} 
\]
where $\delta'_{i} \in \{0,1,2,\ldots, a-1\}$ for $i = j_r,\ldots, j_0 -1$ and $\delta'_{j_r} \neq 0$.
\el

\begin{proof}
The proof is by induction on $r$.   For  $r=1$ we have
\begin{align*}
n_J & = \delta_{j_0} a^{j_0} -  \delta_{j_1} a^{j_1}  \\
& = \left( \delta_{j_0} - 1\right) a^{j_0}  +  \left( a^{j_0}  -   a^{j_1} \right)  -  \left( \delta_{j_1} - 1 \right) a^{j_1} \\
& =  \left( \delta_{j_0} - 1\right) a^{j_0}  + \sum_{i=j_1}^{j_0 -1} (a-1)a^i   -  \left( \delta_{j_1} - 1 \right) a^{j_1} \\
& =  \left( \delta_{j_0} - 1\right) a^{j_0}  + \sum_{i=j_1+1}^{j_0 -1} (a-1)a^i   +  \left( a - \delta_{j_1} \right) a^{j_1}
\end{align*}
and
\[
1 \leq a - \delta_{j_1} \leq a-1.
\]
The nonzero digits of $n_J$ form a single $a$-adic block, and so $M_A(n_J) = 1$.  

Let $r \geq 1$ and suppose that the Lemma is true for $r$.  
Let $J = \{ j_i\}_{i=0}^{r+1}$ be a strictly decreasing sequence  of nonnegative integers and let $\delta_{j_i} \in \{1,2,\ldots, a-1\}$ for $i = 0,1,\ldots, r,r+1$.  We consider the integer 
\[
n_J = \delta_{j_0} a^{j_0} - \sum_{i=1}^{r+1} \delta_{j_i} a^{j_i}.
\]
Let
\[
m_j  = \delta_{j_0} a^{j_0} - \sum_{i=1}^{r} \delta_{j_i} a^{j_i}.
\]
By the induction hypothesis, $M_A(m_J) \leq r$, and $a^{j_r}$ is the smallest power of $a$ that appears in the $a$-adic expansion of $m_J$ with a nonzero digit.  It follows that 
\[
M_A(m_J) - 1 \leq M_A(m_J - a^{j_r}) \leq M_A(m_J) \leq r.
\]
We write
\[
n_J = m_J - \delta_{j_{r+1}} a^{j_{r+1}} 
= \left( m_J - a^{j_r} \right) + \left( a^{j_r} - \delta_{j_{r+1}} a^{j_{r+1}}  \right).
\]
Again applying the induction hypothesis, we see that the positive integer $a^{j_r} - \delta_{j_{r+1}} a^{j_{r+1}}$ has exactly one maximal $a$-adic block, and that the largest power of $a$ that appears in its $a$-adic expansion with a nonzero digit is less than $a^{j_r}$.  
It follows that  
\[
 M_A(m_J) \leq M_A(n_J) \leq M_A(m_J)+1 \leq r+1.
\] 
Moreover, $a^{j_{r+1}}$ is the smallest power of $a$ that appears in the $a$-adic expansion of $a^{j_r} - \delta_{j_{r+1}} a^{j_{r+1}}$ with a nonzero digit,  and so $a^{j_{r+1}}$ is the smallest power of $a$ that appears in the $a$-adic expansion of $n_J$ with a nonzero digit. 
This completes the proof.  
\end{proof}

\bl  \label{quasi:lemma:interpolate}
Let $I$ and $W$ be disjoint finite sets of nonnegative integers.  Let $a \geq 2$, and let $\delta_{i} \in \{1,2,\ldots, a-1\}$ for $i \in I$ and $\delta_w \in \{1,2,\ldots, a-1\}$ for $w \in W$.  Then 
\[
M_A\left( \sum_{i\in I} \delta_i a^i \right) - |W| \leq M_A\left( \sum_{i\in I} \delta_i a^i +  \sum_{w \in W} \delta_w a^w\right) \leq 
M_A\left( \sum_{i\in I} \delta_i a^i \right) + |W|.
\]
\el

\begin{proof}
It suffices to prove the Lemma for $|W|=1$.  
Adding a ``new'' power of $a$ to an $a$-adic representation changes  a  zero digit to a nonzero digit.  If the former zero digit was adjacent to two nonzero digits, then the number of maximal $a$-adic blocks decreases by 1.   If the former zero digit was adjacent to one zero digit and to one nonzero digit, then the number of maximal $a$-adic blocks does not change.   If the former zero digit was adjacent to two zero digits, then the number of maximal $a$-adic blocks increases by 1.  This completes the proof.  
\end{proof}

\bl    \label{quasi:lemma:BoundedBlock}
Let $a \geq 2$ and $k \geq 1$.  If $n$ is a positive integer such that 
\[
n = \sum_{t \in T}\varepsilon_t \delta_t a^t
\]
where $T$ is a  set of $k$ nonnegative integers,  and $\delta_t \in \{1,2,\ldots, a-1\}$ and $\varepsilon_t \in \{1,-1\}$ for all $t \in T$, then $M_A(n) \leq k$.
\el

\begin{proof}
Since $n$ is positive, it follows that $\varepsilon_{t^*}=1$ for $t^* = \max(T)$. 
If $\varepsilon_t = 1$ for all $t \in T$, then $n = \sum_{t \in T}\delta_t a^t$ is the $a$-adic representation, which has exactly $k$ nonzero digits, and so $M_A(n) \leq k$.

Suppose that $\varepsilon_t= -1$ for some $t \in T$.    Arrange $T$ in strictly increasing order $t_1 < t_2 < \cdots < t_k$.  Let $U$ be the set of all $t_i$ with $i \geq 2$ such that  $\varepsilon_{t_i} = 1$ and $\varepsilon_{t_{i-1}} = -1$.  Let $\ell = \card_A(U)$.  We observe that $\varepsilon_{t^*}=1$ implies that $\ell \geq 1$.  Arrange the elements of $U$ in strictly increasing order 
\[
u_1 < \cdots < u_{\ell}.
\]
Define $u_0 = -1$.  For $j = 1, 2,\ldots, \ell$, let
\[
V_j = \{ t\in T : u_{j-1} < t < u_j \text{ and } \varepsilon_t = -1\}.
\]
The set $V_j$ is nonempty for all $j = 1, 2,\ldots, \ell$.  If
\[
v_j = \min(V_j)
\]
then 
\[
u_{j-1} < v_j < u_j  \qquad\text{for $j = 1,\ldots, \ell$.}
\]
Moreover, $\varepsilon_{t^*} =  1$ implies that if $t \in T$ and $\varepsilon_t = -1$, then $t \in V_j$ for some $j$.  
Let $V = \cup_{j=1}^{\ell} V_j$.  We define
\[
n_{V_j} = \varepsilon_{u_j} \delta_{u_j}a^{u_j} + \sum_{v\in V_j}  \varepsilon_v \delta_{v}a^v 
= \delta_{u_j}a^{u_j} - \sum_{v\in V_j} \delta_{v}a^v.
\]
By Lemma~\ref{quasi:lemma:binary}, 
\beq   \label{quasi:UVineq}
a^{u_{j-1}+1} \leq a^{v_j} \leq n_{V_j} < a^{u_j +1}
\eeq
and $a^{v_j}$ is the smallest power of $a$ that appears in the $a$-adic expansion of $n_{V_j}$ with a nonzero digit.  Also,  
\[
M_A\left( n_{V_j} \right) \leq \card_A(V_j).
\]
We define
\[
n' = \sum_{t\in U \cup V} \varepsilon_t \delta_t a^t 
= \sum_{j=1}^{\ell} n_{V_j}.
\]
Then 
\[
M_A(n') \leq  \sum_{j=1}^{\ell} M_A\left( n_{V_j} \right) \leq \sum_{j=1}^{\ell} \card_A(V_j) = \card_A(V).
\]

Let $W = T\setminus (U \cup V).$   
Then $\varepsilon_w = 1$ for all $w \in W$, and 
\[
n = n' + \sum_{w\in W} \delta_w a^w.
\]
Let $I$ be the set of all nonnegative integers $i$ such that $a^i$ occurs with a nonzero digit in the $a$-adic representation of $n'$.  If $i \in I$, then $v_j \leq i \leq u_j$ for some $j  \in \{1,2,\ldots, \ell\}$.  On the other hand, if $w \in W$, then $w > u_{\ell}$ or $u_{j-1} < w < v_j$ for some $j  \in \{1,2,\ldots, \ell\}$.   Therefore, $I \cap W = \emptyset$.   An application of  Lemma~\ref{quasi:lemma:interpolate} gives 
\[
M_A(n) \leq M_A(n') + |W| \leq |V| + |W| = k - |U| \leq	 k-1.
\]
This completes the proof.
\end{proof}

Let $n = \sum_{i=0}^{\infty} \delta_i a^i$ be the $a$-adic expansion of the positive integer $n$.  We introduce the function
\[
\ord_A(n)  = \max\{i : \delta_i \neq 0\}.
\]
Let $r = \ord_A(n)$.  For every positive integer $k \leq r+1$, we call the $k$-tuple 
\[
(\delta_{r-k+1}, \delta_{r-k+2}, \ldots, \delta_{r-1}, \delta_r ) \in \{ 0,1,2,\ldots, a-1\}^{k-1} \times \{ 1,2,\ldots, a-1\}
\]
the \emph{leading $k$-digit string of $n$ with respect to $a$}.

The following result is presumably well known.  

\bl      \label{quasi:lemma:LeadingDigits}
Let $a$ and $b$ be integers greater than 1 such that $a^m \neq b^n$ for all positive integers $m$ and $n$.   
Let $(\gamma_0, \gamma_1, \ldots, \gamma_{k-2}, \gamma_{k-1})$ be a $k$-tuple in $\{ 0,1,2,\ldots, a-1\}^{k-1} \times \{ 1,2,\ldots, a-1\}$.  There exist infinitely many positive integers $n$ such that $b^n$ has leading $k$-digit string  $(\gamma_0, \gamma_1, \ldots, \gamma_{k-2}, \gamma_{k-1})$ with respect to $a$.
\el

\begin{proof}
We claim that the positive real number $\log b/(k\log a)$ is irrational for all positive integers $k$.  If not, then there exist positive integers $r$ and $s$ such that $\log b/(k\log a) = r/s$, or equivalently, $a^{kr} = b^s$, which is absurd.

Let $t = \sum_{i=0}^{k-1} \gamma_i a^i$.   Since $\gamma_{k-1} \in  \{ 1,2,\ldots, a-1\}$, we have $a^{k-1} \leq t < a^k$, and so 
\[
0 \leq \frac{k-1}{k} \leq \frac{\log t}{k\log a} <  \frac{\log (t+1)}{k\log a}\leq 1.
\]

Let $\{ x \}$ denote the fractional part of the real number $x$.  Since the sequence of fractional parts of the positive integral multiples of an irrational number is uniformly distributed in the unit interval $[0,1)$, it follows that there exists a set $\mathcal{N}$ of positive integers of positive asymptotic density $\log((t+1)/t)/(k\log a)$ such that, for every $n \in \mathcal{N}$, we have 
\[
\frac{\log t}{k\log a} \leq \left\{   \frac{n \log b}{k \log a}  \right\} < \frac{\log (t+1)}{k\log a}.
\]
Thus, for every $n \in \mathcal{N}$ there is a positive integer $m$ such that 
\[
\frac{\log t}{k\log a} \leq \frac{n \log b}{k \log a}- m < \frac{\log (t+1)}{k\log a}.
\]
This implies that 
\[
ta^{km} \leq b^n < (t+1)a^{km}
\]
and so $\ord_A\left(b^n - ta^{km}\right) \leq km-1$.  It follows that there exist $\delta_i \in \{0,1,\ldots, a-1\}$ for $i=  0,1,\ldots, km-1$ such that 
\[
b^n - ta^{km} = \sum_{i=0}^{km-1} \delta_i a^i.
\]
Therefore,
\begin{align*}
b^n & = \sum_{i=0}^{km-1} \delta_i a^i +  ta^{km} \\
& = \sum_{i=0}^{km-1} \delta_i a^i +   \left( \sum_{i=0}^{k-1} \gamma_i a^i \right) a^{km} \\
&  = \sum_{i=0}^{km-1} \delta_i a^i 
+   \sum_{i=0}^{k-1} \gamma_i a^{km+i}.
\end{align*}
Thus, $b^n$ has  leading $k$-digit string 
 $(\gamma_0, \gamma_1, \ldots, \gamma_{k-2}, \gamma_{k-1})$ with respect to $a$.  
\end{proof}

\bl    \label{quasi:lemma:LeadingPowerDigits}
Let $a$ and $b$ be integers greater than 1 such that $a^m \neq b^n$ for all positive integers $m$ and $n$.    For every positive integer $\ell$ there exist infinitely many positive integers $n$ such that $M_A(b^n) \geq \ell$.
\el

\begin{proof}
By Lemma~\ref{quasi:lemma:LeadingDigits}, there exist infinitely many positive integers $n$ such that $b^n$ has leading $2\ell$-digit string $(0,1,0,1,\ldots, 0,1,0,1)$ with respect to $a$.  For each such $n$ we have $M_A(b^n) \geq \ell$.
\end{proof}

\section{Proof of the Main Theorem}

\bt      \label{quasi:theorem:Main1}
Let $a$ and $b$ be integers greater than 1, and consider the additive  group \Z\ with generating sets  $A = \{a^i\}_{i=0}^{\infty}$ and $B = \{b^j\}_{j=0}^{\infty}$.  Let $d_A$ and $d_B$ be the metrics on \Z\ associated with the generating sets $A$ and $B$, respectively.   If there exist positive integers $m$ and $n$ such that $a^m = b^n$, then these metrics are bi-Lipschitz equivalent.
\et

\begin{proof}
Since $a^m = b^n$, it follows that $a^{qm} = b^{qn}$  and so $\ell_B(a^{qm}) = 1$ for all nonnegative integers $q$.
We define $H_A = 1+\max\{ \ell_B(a^r) : r=0,1,\ldots, m-1\}$.  

Let $i$ be a nonnegative integer.  By the division algorithm, there exist nonnegative integers $q$ and $r$ such that $i = qm+r$ and $0 \leq r \leq m-1$.  By inequality~\eqref{quasi:LengthIneq},
\[
\ell_B(a^i) = \ell_B\left(a^{qm}a^r\right)  \leq \ell_B\left(a^{qm}\right) + \ell_B\left(a^r\right) \leq 1 + \ell_B\left(a^r\right) \leq H_A.
\]
Therefore, $\sup\{\ell_B(a^i) : a^i \in A \} < \infty$.
Similarly, $\sup\{\ell_A(b^j) : b^j \in B \} < \infty$.
Lemma~\ref{quasi:lemma:MetricEquiv} implies that the metrics $d_A$ and $d_B$ are bi-Lipschitz equivalent.  
This completes the proof.
\end{proof}

\bt    \label{quasi:theorem:Main2}
Let $a$ and $b$ be integers greater than 1, and consider the additive  group \Z\ with generating sets  $A = \{a^i\}_{i=0}^{\infty}$ and $B = \{b^j\}_{j=0}^{\infty}$.  Let $d_A$ and $d_B$ be the metrics on \Z\ associated with the generating sets $A$ and $B$, respectively.  If these metrics are bi-Lipschitz equivalent, then there exist positive integers $m$ and $n$ such that $a^m = b^n$.
\et

\begin{proof}
Since $d_A$ and $d_B$ are bi-Lipschitz equivalent metrics on the group $G$, Theorem~\ref{quasi:theorem:MetricEquiv} implies that 
\[
L = \sup\{ \ell_A(b^j) : b^j \in B \} < \infty
\]
and so every generator $b^j \in B$ can be represented as a word with respect to $A$ of length at most $L$.  By Lemma~\ref{quasi:lemma:BoundedBlock}, $M_A(b^j) \leq L$.  If $a^m \neq b^n$ for all positive integers $m$ and $n$, then Lemma~\ref{quasi:lemma:LeadingPowerDigits} implies that there exist infinitely many $n$ such that $b^n$ has leading $2(L+1)$-digit string $(0,1,0,1,\ldots,0,1)$, and for these numbers $b^n$ we have $M_A(b^n) \geq L+1$. This is a contradiction, and so $a^m = b^n$ for some $m$ and $n$.
\end{proof}

\section{Quasi-isometry}
\bprob
Let $a$ and $b$ be integers greater than 1,  and let $d_A$ and $d_B$ be the metrics on \Z\ associated with the generating sets $A = \{a^i\}_{i=0}^{\infty}$ and $B = \{b^j\}_{j=0}^{\infty}$, respectively.  By the Main Theorem, the identity map from \Z\ to \Z\ is a bi-Lipschitz equivalence if and only if there exist positive integers $m$ and $n$ such that $a^m = b^n$.   It is an open problem to determine if there exists \emph{some} map $f:\Z \rightarrow \Z$ that is a bi-Lipschitz equivalence with respect to these metrics, that is, to determine if the metric spaces $(\Z,d_A)$ and $(\Z, d_B)$ are bi-Lipschitz equivalent.
\eprob

Metric spaces $(X,d_X)$ and $(Y,d_Y)$ are called \emph{quasi-isometric} if there are subspaces $X' \subseteq X$ and $Y' \subseteq Y$ and a number $C > 0$ such that (i) the metric spaces $(X',d_X)$ and $(Y',d_Y)$ are bi-Lipschitz equivalent, and (ii) for every $x\in X$ there exists $x' \in X'$ such that $d_X(x,x') < C$, and for every $y \in Y$  there exists $y'\in Y'$ such that $d_Y(y,y') < C$.  

\bprob
Let $a$ and $b$ be integers greater than 1, and let $d_A$ and $d_B$ be the metrics on \Z\ associated with the generating sets $A = \{a^i\}_{i=0}^{\infty}$ and $B = \{b^j\}_{j=0}^{\infty}$, respectively.  Are the metric spaces $(\Z,d_A)$ and $(\Z, d_B)$ quasi-isometric?  This problem was first posed by Richard E. Schwartz for the generating sets $A = \{ 2^i \}_{i=0}^{\infty}$ and $B  = \{ 3^j \}_{j=0}^{\infty}$.  
\eprob

\emph{Acknowledgements.}  I am grateful to Jason Behrstock for introducing me to Schwartz's question about arithmetic quasi-isometry, which led directly to the problem of the bi-Lipschitz equivalence of geometric generating sets for \Z\ that is solved in this paper.   I also thank Jacob Fox for many helpful contributions to this work.

\end{document}